\documentclass{amsart}

\usepackage{amsmath,amssymb}

\newtheorem{lemma}{Lemma}
\newtheorem{prop}[lemma]{Proposition}
\newtheorem{cor}[lemma]{Corollary}
\newtheorem{thm}[lemma]{Theorem}

\newtheorem{thm?}[lemma]{Theorem?}

\newtheorem*{mainthm}{Main Theorem}

\newtheorem*{unthm}{Theorem}
\newtheorem*{fact}{Fact}

\newcommand{\ra}{\ensuremath{\rightarrow}}

\newcommand{\Z}{\mathbb{Z}}
\newcommand{\N}{\mathcal{N}}
\newcommand{\R}{\mathbb{R}}

\newcommand{\C}{\mathbb{C}}
\newcommand{\Q}{\mathbb{Q}}

\newcommand{\ab}{\operatorname{ab}}

\newcommand{\End}{\operatorname{End}}

\newcommand{\Gal}{\operatorname{Gal}}
\newcommand{\GL}{\operatorname{GL}}

\newcommand{\tors}{\operatorname{tors}}

\newcommand{\Pic}{\operatorname{Pic}}
\newcommand{\Frac}{\operatorname{Frac}}

\newcommand{\Spec}{\operatorname{Spec}}
\newcommand{\Prin}{\operatorname{Prin}}

\newcommand{\car}{\operatorname{char}}
\newcommand{\Div}{\operatorname{Div}}
\newcommand{\Cl}{\mathbf{Cl}}
\newcommand{\MW}{\bf{MW}}
\newcommand{\MWQ}{\bf{MWQ}}
\newcommand{\MWT}{\bf{MW/T}}
\newcommand{\tf}{\operatorname{tf}}
\newcommand{\sep}{\operatorname{sep}}
\renewcommand{\car}{\operatorname{char}}

\address{Department of Mathematics \\ Boyd Graduate Studies Research Center \\ University 
of Georgia \\ Athens, GA 30602-7403 \\ USA}
\email{pete@math.uga.edu}

\title{Elliptic Dedekind Domains Revisited}

\author{Pete L. Clark}

\begin{document}
\maketitle

\begin{abstract}
We give an affirmative answer to a 1976 question of M. Rosen: every abelian group is isomorphic to the class group of 
an elliptic Dedekind domain $R$.  We can choose $R$ to be the integral closure of a 
PID in a separable quadratic field extension.  In particular, this yields new and -- we feel -- simpler proofs of 
theorems of L. Claborn and C.R. Leedham-Green.
\end{abstract} \noindent \\ \\
Luther Claborn received his PhD from U. Michigan in 1963 and died 
in a car accident in August of 1967.  In between he wrote $11$ papers, including \cite{Claborn3}, which shows that every abelian group is the class group of a Dedekind domain.  Evidently this result remains of interest 
to this day, more than $40$ years after his untimely death.  This paper is dedicated to him.
\\ \\
Terminology: For a scheme $X$, $\Pic X$ denotes the Picard group $H^1(X,\mathcal{O}_X^{\times})$, whose elements 
correspond to isomorphism classes of line bundles on $X$.  When $X = \Spec R$ is affine, we abbreviate $\Pic(\Spec R)$ to $\Pic(R)$.  When 
$R$ is a Dedekind domain, $\Pic(R)$ is the \textbf{ideal class group} of $R$. \\ \indent
An \textbf{overring} of an integral domain $R$ is a ring intermediate between $R$ and its fraction field.  A Dedekind 
domain $R$ will be called \textbf{affine} if it is the coordinate ring of a nonsingular, geometrically integral 
affine curve $C$ defined over some field $k$.  A \textbf{geometric} Dedekind domain is an overring of an affine domain.  \\ 
\indent
An elliptic curve over a field $k$ means, as usual, a complete, nonsingular, geometrically integral 
genus one curve $E_{/k}$ with a distinguished $k$-rational point $O$.  To an elliptic curve we associate its 
function field $k(E)$ and its \textbf{standard affine ring} $k[E]$, the ring of all functions on $E$ which are regular 
away from $O$.  An \textbf{elliptic domain} is an overring of the standard affine ring of some elliptic curve.  \\ \indent
A \textbf{quadratic domain} is a Dedekind domain which can be obtained by taking the integral closure of a PID in a 
quadratic field extension. \\ \indent
For a cardinal $\kappa$, $FA(\kappa) := \bigoplus_{x \in \kappa} \Z$ denotes the free abelian group of rank $\kappa$.

\section{Introduction}
A celebrated 1966 theorem of Luther Claborn asserts that for any abelian group $A$ whatsoever, there exists a Dedekind 
domain $R$ whose ideal class group $\Pic(R)$ is isomorphic to $A$.  A different proof was given in 1972 by C.R. 
Leedham-Green, which shows that $R$ may be taken to be quadratic, i.e., the integral closure of a PID in a 
quadratic field extension.  Claborn's proof makes use of some nontrivial facts from commutative algebra; especially, it requires 
familiarity with the divisor class group of a Krull domain.  Leedham-Green's proof is more elementary -- in his own words, it is ``based on a naive geometrical construction'' -- but 
is quite intricate.  
\\ \\ \noindent
Work of M. Rosen takes a completely different approach, based on the following:
\\ \\
Fact I. For an elliptic curve $E$ over a field $k$, the standard affine ring $R = k[E]$ is a Dedekind domain with 
$\Pic(R)$ isomorphic to the Mordell-Weil group $E(k)$.
\\ \\
Fact II. With $R = k[E]$ as above, for any subgroup $H \subset \Pic(R)$, there exists an overring $R^H$ of $R$ such that $\Pic(R^H) \cong \Pic(R)/H$.  
\\ \\
From these two facts it follows that any abelian group which is isomorphic to a quotient group of a Mordell-Weil group 
is the class group of some Dedekind domain.  Rosen calls a Dedekind domain which arising as an overring of the standard affine 
ring of some elliptic curve \textbf{elliptic}.  
\\ \\
In \cite{Rosen73}, Rosen shows that any finitely generated abelian group is the class group of the coordinate of 
some (not necessarily standard) affine elliptic curve over some number field $k$.  In \cite{Rosen76}, Rosen uses Serre's open image theorems to show that 
every countably generated abelian group is the class group of an elliptic Dedekind domain.  His 
method does not work for uncountable groups, and accordingly he asks whether every abelian group is the class group of 
an elliptic Dedekind domain.
\\ \\
Our main result gives an affirmative answer to this question.
\begin{mainthm} \textbf{} \\ a) For an abelian group $A$, there is an elliptic Dedekind domain $R$ with $\Pic(R) \cong A$. \\
b) We can moreover choose the elliptic Dedekind domain $R$ to be quadratic, i.e., the integral closure of a PID in a 
quadratic field extension.
\end{mainthm}
\noindent
Our construction is most certainly inspired by Rosen's work, and it follows his general strategy in that it uses Facts I. and II. above to reduce to the problem of  
constructing a free abelian group of arbitrary rank as a quotient of some Mordell-Weil group.  But there are also 
several differences.  First, we in fact construct arbitrary free abelian groups as Mordell-Weil groups, whereas 
Rosen constructs a free abelian group of countable rank as the quotient of a Mordell-Weil group by its torsion subgroup.  
Second, whereas Rosen's construction takes $k$ to be the maximal multiquadratic extension of $\Q$, ours does not.  Nor could it, of 
course: the group of rational points of an elliptic curve over a countable field must be countable.  Our field $k$ is a (transfinitely!) iterated 
function field, and accordingly we make no use of Serre's open image theorems nor any other deep arithmetic facts.  
\\ \\
The proof of the Main Theorem occupies little more than a single page.  However, our goal is to give a proof which a student 
or practitioner of arithmetic algebraic geometry will find self-contained and transparent.  To this end, we have included 
in $\S 2$ some material on class groups of overrings of Dedekind domains that is of an expository nature.  We have 
also included a short discussion of the property of a Dedekind domain that every ideal class be represented by at least 
one prime ideal -- we call such a domain \textbf{replete} -- as well as a slightly weaker property that suffices for 
the intended applications.  Our Theorem \ref{2.12} on the repleteness and weak repleteness of elliptic Dedekind domains 
may be new.  We give the proof of the Main Theorem in $\S 3$.

\section{Preliminaries}
\noindent
We wish to recall some results concerning the effect of passage to an overring on the class 
group, and on the connection between class groups of affine curves and Picard groups of the Jacobians 
of their projective completions.  We could not resist mentioning a few interesting results which are closely related to these 
topics but not needed for the proof of the Main Theorem.  For such results we explicitly state that they are not 
needed in the sequel, and we give references rather than proofs.
\subsection{Basic definitions} \textbf{} \\ \\
\noindent For a Dedekind domain $R$, let $\Sigma(R)$ be the set of all nonzero prime ideals; we typically speak of 
elements of $\Sigma(R)$ as simply ``primes.'' Consider the map 
\[\Phi: \Sigma \ra \Pic(R), \ \mathfrak{p} \mapsto [\mathfrak{p}]. \]
Since the group $\Frac(R)$ of fractional ideals of $R$ is free abelian with $\Sigma$ as a basis, 
$\Phi$ uniquely extends to a homomorphism $\Frac(R) \ra \Pic(R)$, which is surjective, and whose kernel is the subgroup 
$\Prin(R)$ of principal fractional ideals.  
\\ \\
If $R$ and $S$ are Dedekind domains, by a \textbf{morphism} of Dedekind domains we mean an injective ring homomorphism 
$\iota: R \hookrightarrow S$.  If $I$ is a fractional ideal of $R$, then the \textbf{push forward} $I \in \Frac(R) 
\mapsto IS$ induces a homomorphism from $\Frac(R)$ to $\Frac(S)$, denoted $\iota_*$.  Since the pushforward of a principal 
fractional ideal remains principal, $\iota_*$ factors through to a homomorphism $\iota_*: \Pic(R) \ra \Pic(S)$.

\subsection{Overrings} \textbf{} \\ \\
\noindent
%
If $R$ is an integral domain with field of fractions $K$, an \textbf{overring} of $R$ is a ring $S$ intermediate between 
$R$ and $K$, i.e., $R \subset S \subset K$.

\begin{lemma}
\label{2.1}
Let $\iota: R \hookrightarrow S$, where $R$ is a Dedekind domain and $S$ is an overring. \\
a) For any $\mathcal{P} \in \Sigma(S)$, $S_{\mathcal{P}} = R_{\mathcal{P} \cap R}$.  \\
b) $S$ is itself a Dedekind domain. \\
c) $\iota^*: \Sigma(S) \hookrightarrow \Sigma(R)$. \\
d) For all $\mathcal{P} \in \Sigma(S)$, $\iota_* (\iota^* \mathcal{P}) = \mathcal{P}$. 
\end{lemma}
\noindent
Proof: a) Put $\mathfrak{p} = \mathcal{P} \cap R$.  If $0 \neq \frac{x}{y} \in \mathcal{P}$, then $0 \neq x = y (\frac{x}{y}) \in \mathfrak{p}$, so 
$\mathfrak{p} \in \Sigma(R)$.  Thus $S_{\mathcal{P}}$ contains the DVR $R_{\mathfrak{p}}$ and is properly contained in its 
fraction field, so $S_{\mathcal{P}} = R_{\mathfrak{p}}$. \\ \indent b)  By the Krull-Akizuki Theorem \cite[Thm 11.7]{Matsumura}, $S$ 
is a one-dimensional Noetherian domain; and by part a) the localization of $S$ at every prime is a DVR, hence $S$ is integrally closed and 
thus a Dedekind domain.  \\ \indent c) From part a), we have that 
there is no other prime $\mathcal{P}'$ of $S$ with $\iota^*(\mathcal{P}') = \mathfrak{p}$, since localizations at 
distinct primes in a Dedekind ring are distinct DVRs. \\ \indent d) Moreover 
 \[\iota_*(\iota^*(\mathcal{P})S_{\mathcal{P}})) = S_{\mathcal{P}} \mathfrak{p} R = \mathcal{P} S_{\mathcal{P}}. \]
By part b), $\iota_*(\mathfrak{p} R)$ is not divisible by any prime other than $\mathcal{P}$, so that 
$\iota_*(\mathfrak{p}) = \mathcal{P}$. 
\begin{cor}
\label{2.2}
Let $S$ be an overring of the Dedekind domain $R$.  The prime ideals of $S$ are identified, via $\iota^*$, with 
the prime ideals $\mathfrak{p}$ of $R$ such that $\mathfrak{p}S \subsetneq S$.
\end{cor}
\noindent
We can explicitly describe all overrings of a Dedekind domain $R$.  For an arbitrary subset $W \subset \Sigma(R)$, 
put $R_W := \bigcap_{\mathfrak{p} \in W} R_{\mathfrak{p}}$, the intersection taking place in the fraction field $K$.  
For the sake of simplifying some later formulas, we also define 
\[R^W = R_{\Sigma(R) \setminus W}. \]
Evidently $R_W$ is an overring of $R$, so is itself a Dedekind domain.  Conversely:
\begin{thm}
\label{2.3}
Let $R$ be a Dedekind with fraction field $K$, and let $R \subset S \subset K$ be an overring.  Let 
$W$ be the set of all primes $\mathfrak{p}$ of $R$ such that $\mathfrak{p}S \subsetneq S$. Then \[S = 
R_W = \bigcap_{\mathfrak{p} \in W} R_{\mathfrak{p}}. \]
\end{thm}
\noindent
Proof: Not needed in the sequel; see e.g. \cite[Cor. 6.12]{LM}.
\\ \\
The reader may be more used to thinking about generating overrings by \emph{localization}: if $R$ is a domain 
and $T \subset R$ a multiplicatively closed set, then $R[T^{-1}]$ is an overring of $R$.  Probably the reader knows 
that every overring of $\Z$ is obtained by localization; in fact this is true for overrings of any PID $R$.  For 
this it suffices to exhibit $R[\frac{x}{y}]$ as $R[\frac{1}{z}]$.  The key point here is that since $R$ is a UFD 
we can assume that $x$ and $y$ are relatively prime, and then there exist $a,b \in R$ with $ax + by = 1$, so 
that $\frac{1}{y} = \frac{ax+by}{y} = a \left(\frac{x}{y} \right) + b \left( \frac{y}{y} \right) \in R[\frac{x}{y}]$.\footnote{The 
fact that there is no ring strictly intermediate between a DVR and its fraction field, which was used in the proof of 
Lemma \ref{2.1}, is an (even) easier special case.}
But in general, not all overrings are realizable by localization:
\begin{thm} 
\label{2.4}
Let $R$ be a Noetherian domain, and consider the following properties: \\
(i) Every overring of $R$ is integrally closed. \\
(ii) Every overring of $R$ is obtained by localizing at a multiplicative subset. \\
Then (i) holds iff $R$ is a Dedekind domain, and (ii) holds iff $R$ is a Dedekind domain with torsion class group.
\end{thm}
\noindent
Proof: Not needed in the sequel; see \cite{Davis} or \cite{GO}.
\\ \\
The following result explains the importance of overrings in the study of class groups of Dedekind domains.
\begin{thm}(Claborn, \cite{Claborn2})
\label{2.5}
Let $R$ be a Dedekind domain, and $S = R^W$ be an overring of $R$.  There exists a short exact sequence 
\[0 \ra H \ra \Pic(R) \stackrel{\iota_*}{\ra} \Pic(S) \ra 0, \]
where $H = \langle \Phi(W) \rangle$ is the subgroup generated by classes of primes $\mathfrak{p}$ 
with $\mathfrak{p} S = S$.
\end{thm}
\noindent
Proof: Since $\iota_* \circ \iota^* = 1_{\Sigma(S)}$, $\iota_*$ is surjective on prime ideals; \emph{a fortiori} 
the induced map on class groups is surjective.  Clearly each prime $\mathfrak{p}$ with $\mathfrak{p}S = S$ lies 
in the kernel.  Conversely, suppose $I$ is a fractional ideal of $R$ in the kernel of $\iota_*$, so that 
there exists $x$ in the fraction field with $IS = xS$.  Then $xI^{-1}S = S$, so that $xI^{-1}$ is a product of 
primes $\mathfrak{p}$ with $\mathfrak{p}S = S$.
\\ \\
Thus one can realize certain quotients of the the Picard group of $R$ by passing to a suitable overring $S$.  In general 
however, not every subgroup of $\Pic(R)$ is generated by classes of prime ideals.  This brings us to the next section.

\subsection{Replete Dedekind rings} \textbf{} \\ \\
\noindent
We say that a Dedekind domain $R$ is \textbf{replete} if the map $\Phi$ is surjective, i.e., if every ideal 
class is represented by a prime ideal.
\begin{prop}
\label{2.6}
Let $R$ be a replete domain, and let $H \subset \Pic(R)$ be any subgroup.  Then there exists an overring 
$S$ of $R$ such that $\Pic(S) \cong \Pic(R)/H$.
\end{prop}
\noindent
Proof: Indeed, if $R$ is replete, then $H$ is generated by a set $W$ of prime ideals of $R$.  Then take 
$S$ to be the overring $R^W = \bigcap_{\mathfrak{p} \in \Sigma(R) \setminus W} R_{\mathfrak{p}}$.  By Theorem \ref{2.5}, 
$\Pic(R^W) \cong \Pic(R)/H$.
\\ \\
For the proof of Proposition \ref{2.6} to go through, it suffices that $R$ have the property that any subgroup $H$ of 
$\Pic(R)$ is generated by classes of prime ideals.  Let us call a domain with this property \textbf{weakly replete}.
\begin{cor} 
\label{2.7}
An overring of a weakly replete domain is weakly replete.
\end{cor}
\noindent
Proof: This follows easily from Theorem \ref{2.5}.
\\ \\
Examples: Trivially a PID is replete.  The repleteness of the ring of integers in a global field is a weak version of 
the Chebotarev Density Theorem.  We will see in $\S 2.4$ that the standard affine ring of an elliptic curve is weakly 
replete but not necessarily replete.  Examples of domains which are not weakly replete seem harder to come by.  In \cite{Claborn4}, 
Claborn exhibits for each $n \in \Z^+$ a Dedekind domain $R_n$ whose class group is cyclic of order $n$ and such that $[\mathfrak{p}] = 
[\mathfrak{q}]$ for all $\mathfrak{p}, \mathfrak{q} \in \Sigma(R_n)$ as well as a Dedekind domain $R$ with 
$\Pic(R) \cong \Z$ such that for all $\mathfrak{p} \in \Sigma(R)$, $[\mathfrak{p}] = \pm 1$.  
\\ \\
A \textbf{repletion} of a Dedekind domain $R$ is a replete Dedekind domain $S$ together with an injective ring 
homomorphism $\iota: R \hookrightarrow S$, such that $\iota_*: \Pic(R) \stackrel{\sim}{\ra} \Pic(S)$.  
\begin{thm}(Claborn)
\label{2.8}
For a Dedekind domain $R$, let $R^1$ denote the localization of $R[t]$ at the 
multiplicative set generated by all monic polynomials.  Then $R^1$ is Dedekind and the composite map $\iota: R \rightarrow R[t] \rightarrow R^1$ is a 
repletion.
\end{thm}
\noindent
Proof: Not needed in the sequel; see \cite[Cor. 2.5]{Claborn1}.
%
\begin{cor}
\label{2.9}
(Claborn) 
For any Dedekind domain $R$, and any subgroup $H \subset \Pic(R)$, there exists a Dedekind domain $S$ and a 
homomorphism of Dedekind domains $\iota: R \ra S$ making the following sequence exact: 
\[0 \ra H \subset \Pic(R) \stackrel{\iota_*}{\ra} \Pic(S) \ra 0. \]
Thus every quotient group of $\Pic(R)$ is the class group of some Dedekind domain.
\end{cor}
\noindent
Proof: This follows immediately from Theorem \ref{2.8} and Proposition \ref{2.6}.

\subsection{Affine domains, gometric domains, and elliptic domains} \textbf{} \\ \\ \noindent
Let $k$ be a field.  To a pair $(C,O)$, where $C_{/k}$ is a 
complete, nonsingular geometrically integral curve and $O \in C(k)$ is a rational point, we attach the rational 
function field $k(C)$ and \textbf{standard affine ring} $k[C^o]$, the subring of $k(C)$ consisting of all functions 
which are regular on all points except (possibly) $O$.  Note that $k[C]$ is the coordinate ring of the 
affine algebraic curve $C^o = C \setminus O$.  The ring $k[C^o]$ is a nonsingular Noetherian domain of dimension one, i.e., a 
Dedekind domain.  Consider the map which sends a degree $0$ divisor  $\sum_{P} n_P [P]$ on $C$ to the divisor 
$\sum_{P \neq O} n_P [P]$ (of degree $-n_O$) on $C^o$.  Upon quotienting out by principal divisors, this gives an isomorphism
\begin{equation}
\label{GEOISO}
 J(C)(k) = \Pic^0(C) \stackrel{\sim}{\ra} \Pic(k[C^o]),
\end{equation}
where $J$ is the Jacobian of $C$.  Thus the class group of a standard affine domain is canonically isomorphic to 
the group of $k$-rational points on a certain (Jacobian) abelian variety.  
\\ \\
When $C = E$ has genus one, the automorphism group acts transitively on the set of $k$-rational points, so the affine 
curve $E^o$ is independent of the choice of $O$.  In this case, we simplify the notation $k[E^o]$ to $k[E]$.
\\ \\
In general, let us say that a Dedekind domain $R$ is \textbf{affine} if it is of the form $k[C^o]$ for some nonsingular, 
geometrically integral affine curve $C^o$ over a field $k$.  Write $C$ for the nonsingular projective model of 
$C^o$.  As long as $C \setminus C^o$ contains at least one $k$-rational point, a well-known argument using Riemann-Roch shows 
that the affine domain $k[C^o]$ is an overring of some standard affine domain. 
\begin{thm} (Rosen)
\label{2.10} Let $C^o = C \setminus S$ be a nonsingular, geometrically integral affine curve over a field $k$.  Let 
$D^0(S)$ be the subgroup of $\Div(C)$ consisting of degree $0$ divisors supported on $S$, and let $P(S)$ be the principal 
divisors in $D^0(S)$.  Let $d$ be the least positive degree of a divisor supported on $S$ (note that $d = 1$ iff 
$S$ contains at least one $k$-rational point), and let $i$ be the least 
positive degree of a divisor on $C$.  Then there is an exact sequence 
\[0 \ra D^0(S)/P(S) \ra \Pic^0(C) \ra \Pic(C^o) \ra Z(d/i) \ra 0, \]
where $Z(d/i)$ is a cyclic group of order $d/i$.  
\end{thm}
\noindent Proof: Not used in the sequel; see \cite{Rosen73}.
\\ \\
Using this exact sequence and the fact that every elliptic curve over $\overline{\Q}$ has infinite rank, Rosen deduces:
\begin{thm}
\label{2.11}
(Rosen)
For any finitely generated abelian group $A$, there is a number field $k$ and a (not necessarily standard) 
affine elliptic curve $E^o$ over $k$ such that $\Pic(k[E^o]) \cong A$.
\end{thm}
\noindent
The claim that in Theorem \ref{2.11} we can always take $k = \Q$ is equivalent to the existence of elliptic curves $E_{/\Q}$ 
of arbitrarily large rank, a notorious open problem.  
\\ \\
A Dedekind domain is \textbf{geometric} if it is an overring of an affine Dedekind domain.  
In other words, a geometric Dedekind domain is the ring of all functions on a complete curve $C_{/k}$ which are regular 
on some fixed, but possibly infinite, subset of closed points of $C$.  Finally, an \textbf{elliptic} Dedekind domain 
is an overring of the standard affine domain of an elliptic curve $E_{/k}$.  
\begin{thm} 
\label{2.12}
Let $E_{/k}$ be an elliptic curve with equation $y^2 = P(x) = x^3+Ax+B$.  
a) The standard affine ring $k[E]$ is weakly replete (hence so are all of its overrings). \\
b) If $k$ is algebraically closed, $k[E]$ is not replete. \\
c) Suppose $k$ does not have characteristic $2$ and $k[E]$ is not replete.  Then for all $x \in k$, there 
exists $y \in k$ with $y^2 = P(x)$.
\end{thm} 
\noindent
Proof: Each point $P \neq 0$ on $E(k)$ a prime ideal in the standard affine ring $k[E]$; according to the isomorphism of 
(\ref{GEOISO}), every nontrivial element of $\Pic(k[E])$ arises in this way.  This proves part a).  Part b) is similar: 
if $k$ is algebraically closed, then by Riemann-Roch every prime ideal of $k[E]$ corresponds to a $k$-valued point 
$P \neq O$ on $E(k)$, which under (\ref{GEOISO}) corresponds to a nontrivial element of the class group.  Therefore 
the trivial class is not represnted by any prime ideal.  Under the hypotheses of part c), there exists an $x \in k$ 
such that the points $(x, \pm \sqrt{P(x)})$ form a Galois conjugate pair.  Therefore the divisor 
$(x,\sqrt{P(x)})+ (x,-\sqrt{P(x)})$ represents a closed point on the curve $C^o$, in other words a nonzero prime ideal 
of $k[E]$.  But the corresponding point on $E(k)$ is $(x,\sqrt{P(x)}) + (x,-\sqrt{P(x)}) = O$.
\\ \\
To sum up: since every abelian group $A$ is a quotient of a free abelian group $FA(\kappa)$ of some rank $\kappa$, 
and the standard affine domain $k[E]$ attached to an elliptic curve $E_{/k}$ is weakly replete, in order to 
realize $A$ as the Picard group of an elliptic Dedekind domain it suffices to find $k$ and $E$ such that 
$E(k) \cong FA(\kappa)$.  This we handle in the next section, along with the claim that the domain can be taken to be 
the integral closure of a PID in a quadratic extension.

\section{Proof of the Main Theorem}
\noindent
\begin{prop}
\label{SILVERPROP}
Let $K$ be a field and $E_{/K}$ an elliptic curve.  Let $K(E)$ be the function field of $E$.  Then there is a short exact 
sequence 
\[0 \ra E(K) \ra E(K(E)) \ra \End_K(E) \ra 0. \]
Here $\End_K(E) \cong \Z^{a(E)}$, where $a(E) = 2$ if $E$ has $K$-rational CM, and otherwise $a(E) = 1$.  Since $\End_K(E)$ is free abelian, we have 
$E(K(E)) \cong E(K) \bigoplus \Z^{a(E)}.$
\end{prop}
\noindent
Proof: $E(K(E))$ is the group of rational maps from the nonsingular curve $E$ to the complete variety $E$ (the group law 
is pointwise addition).  But every rational map from a nonsingular curve to a complete variety is everywhere defined, so 
$E(K(E))$ is the group of all morphisms $E \ra E$ under pointwise addition.  The constant morphisms form a subgroup 
isomorphic to $E(k)$, and every map of curves from $E$ to itself differs by a unique constant from 
a map of elliptic curves $(E,O) \ra (E,O)$, i.e., an endomorphism of $E$.  
\\ \\
Now take $(E_0)_{/\Q}: y^2+y = x^3-49x - 86$, so $E_0(\Q) = 0$ \cite[Theorem H]{Kolyvagin}.  This elliptic curve has nonintegral $j$-invariant $\frac{2^{12} 3^3}{37}$, so 
does not have complex multiplication.   So defining $K_0 = \Q$ and 
$K_{n+1} = K_n(E_{/K_n})$, Proposition \ref{SILVERPROP} 
gives 
\[ E(K_n) \cong \bigoplus_{i=1}^n \Z. \]
\\ Now define $K_{\infty} = \lim_{n \ra \infty} K_n$; what can we say 
about $E(K_{\infty})$?  We have the following technical result:
\begin{lemma}(``Continuity Lemma'')
Let $K$ be a field, $(K_i)_{i \in I}$ be a directed system of field extensions of $K$, and $E_{/K}$ and elliptic curve.  
Then there is a canonical isomorphism 
\[ \lim_I E(K_i) = E(\lim_I K_i). \]
\end{lemma}
\noindent
Proof: E.g. by general nonsense: this holds for any representable contravariant functor from the category of affine $K$-schemes to the category of abelian groups.  
\\ \\
Therefore $E(K_{\infty}) = \lim_n E(K_n) = \bigoplus_{n \in \Z^+} \Z$, recovering Rosen's Theorem.
\\ \\
Now given an uncountable $S$, choose $\omega$ an ordinal of the same cardinality.  We define the field $K_{\omega}$ by transfinite induction: 
$K_0 = \Q$, for an ordinal $o < \omega$, $K_{o+1} = K_o(E_{/K_0})$, and for a limit ordinal $o$, $K_o = \lim_{o' < o} 
K_{o'}$.  By the Continuity Lemma, we have $E(K_o) = \lim_{o' \in o} E(K_{o'})$.  \\ \indent An isomorphism from 
$E(K_o)$ to $\bigoplus_{o' \in o} \Z$ can be built up by transfinite induction as well; this amounts to the following 
elementary exercise (c.f. \cite[p. 105]{Scott}): 
\begin{fact}
For an abelian group $A$, TFAE: \\
(i) $A$ is free abelian. \\
(ii) $A$ has a well-ordered ascending series with all factors $A_{s+1}/A_s$ infinite cyclic.
\end{fact}
\noindent
Thus for a given abelian group $A = \Z[\kappa]/H$, we have constructed a field $k$, an elliptic curve $E_{/k}$, and an 
overring $R$ of the affine domain $k[E]$ such that $\Pic(R) \cong \Z[\kappa]/H \cong A$, which proves part a) of 
the Main Theorem.
\\ \\
As for the second part, let $\sigma$ be the automorphism of the function field $k(E)$ induced by $(x,y) \mapsto (x,-y)$, 
and notice that $\sigma$ corresponds to inversion $P \mapsto -P$ on $E(k) = \Pic(k[E])$.  Let $S = R^{\sigma}$ be 
the subring of $R$ consisting of all functions which are fixed by $\sigma$.  Then $k[E]^{\sigma} = k[x]$ is a PID, and $S$ is an overring of $k[x]$, hence also a PID.  More precisely, $S$ is the 
overring of all functions on the projective line which are regular away from the point at infinity and the $x$-coordinates 
of all the elements in $H$ (note that since $H$ is a subgroup, it is stable under inversion).  Finally, to see 
that $R$ is the integral closure of $S$ in the quadratic field extension $k(E)/k(x)$, it suffices to observe the following:
For any separable quadratic field extension $L/K$ with nontrivial automorphism $\sigma$, and $R$ a Dedekind domain 
with fraction field $L$, putting $S = R \cap K$ we have that $R$ is the integral closure of $S$ in $L$, since each element 
in $x \in R$ satisfies the polynomial $(t-x)(t-\sigma(x))$.

\end{document}